\documentclass[11pt]{article}

\usepackage[a4paper,margin=0.75in]{geometry}
\usepackage{amsmath,amssymb,amsfonts}
\usepackage{graphicx}
\usepackage{bm}
\usepackage{siunitx}
\usepackage{hyperref}
\usepackage{authblk}
\usepackage{mathtools}
\usepackage{caption}
\usepackage{float}
\usepackage{booktabs}
\usepackage{setspace}
\usepackage{array}
\usepackage{xcolor}
\hypersetup{hidelinks,
  pdftitle={From the Deterministic Hook Effect to Stochastic Reliability in Ternary Complex Formation},
  pdfauthor={Deepak K. Agrawal},
  pdfsubject={Stochastic reliability in ternary complex formation}}
\graphicspath{{figures/}}
\title{From the Deterministic Hook Effect to Stochastic Reliability in Ternary Complex Formation}

\author{Deepak K. Agrawal}
\affil{Department of Biosciences and Bioengineering, Indian Institute of Technology Bombay, Mumbai, Maharashtra, India\\
\texttt{dkagrawal@iitb.ac.in}}

\date{}

\begin{document}
\maketitle

\begin{abstract}
Ternary binding systems occur in signaling, scaffold-mediated assembly, and targeted protein degradation. A mediator brings two partners into a productive ternary complex, but excess mediator can hold them in separate binary complexes. Consequently, ternary complex formation can increase and then decrease as the mediator dose increases, producing the hook effect. Deterministic dose-response curves describe the average behavior but do not indicate whether individual small compartments reliably form the complex. A minimal reversible binding network was studied using deterministic equations, exact stochastic simulations, and a simple catalytic extension. The deterministic model recovered the hook-shaped response. In stochastic simulations, small volumes often failed to accumulate the required number of complexes, even where the deterministic prediction was favorable. At low mediator doses, failure reflected a shortage of mediator molecules; at high doses, material accumulated in free or binary-bound forms. Larger volumes increased the range of mediator doses supporting reliable formation, while stricter copy-number thresholds narrowed that range. Faster ternary-complex completion or slower breakup improved reliability. In the catalytic model, product formation also depended on how early productive complexes formed. These results show why a favorable average binding response does not necessarily ensure reliable assembly in small reaction volumes.
\end{abstract}

\noindent\textbf{Keywords:} ternary binding; hook effect; stochastic simulation; finite-copy reliability; biochemical reaction networks; targeted protein degradation
\section{Introduction}

Ternary complex formation is a common molecular motif in which a mediator binds two partners, bringing them together into a single ternary complex. This motif appears in many systems. Bispecific antibodies and T-cell engagers connect immune and target cells through two binding interactions \cite{klein2012bispecific,zahavi2020monoclonal}. Scaffold proteins bring signaling proteins together and can change pathway timing and specificity \cite{good2011scaffold}. PROTACs and molecular glues recruit a target protein and an E3 ligase into a ternary complex that can drive degradation \cite{chamberlain2019structural,zhao2023kinetic}. Synthetic scaffolds use the same principle to control metabolic flux \cite{dueber2009synthetic}. The same binding logic also underlies sandwich immunoassays and lateral flow assays. At an intermediate analyte concentration, the capture and detection molecules are brought together, increasing signal. When the analyte concentration becomes too high, however, the signal can decrease because the two binding partners are saturated separately rather than forming a productive complex. This non-monotonic response is known as the hook effect or prozone effect \cite{wild2013hook,rey2017mitigating}.

Although the participating molecules differ across these applications, the same competition between productive ternary binding and separate binary binding underlies their dose response. A minimal binding network can therefore help explain the shared constraint without representing every detail of a particular biological system. 

Earlier equilibrium studies showed that the hook-shaped response follows from binding and conservation of material \cite{douglass2013comprehensive,ha2016mathematical,perelson1980equilibrium}. At low mediator doses, there are too few molecules to bring the partners together. At high doses, the partners are increasingly held in separate binary complexes \cite{levchenko2000scaffold,park2023cooperativity,kostic2020critical}. Therefore, productive assembly is strongest over an intermediate dose range. For targeted degradation and other ternary systems, binding affinity, molecular geometry, and interactions within the resulting complex can also influence the outcome \cite{hughes2017protacs,guzzetti2023integrated}. These findings motivate looking beyond the location of the hook-curve peak: a useful mediator dose must also reliably produce the complex.

Most studies describe the hook effect using average concentrations. Such curves show where ternary-complex formation is expected to be high, but not how often an individual small compartment succeeds. This matters in femtoliter-scale volumes and other confined environments, where nanomolar concentrations can correspond to only a few molecules \cite{bhalla2004small,takahashi2005space}. A nonzero predicted mean can coexist with many compartments containing no productive complex, as is familiar in stochastic biochemical systems \cite{hahl2016comparison,elowitz2002stochastic,mcadams1997stochastic,anderson2014stochastic}.

This study links the expected hook-shaped response to the probability of forming productive complexes within a single small compartment. Deterministic calculations define the average response, while exact stochastic simulations show how often a chosen copy-number threshold is reached and whether it is reached within a specified time window. The analysis then tests whether reliable assembly improves by speeding up ternary completion, slowing complex breakup, or favoring one binary binding route. Finally, a catalytic extension examines whether the assembled complex can generate product during a fixed observation period \cite{bondeson2015catalytic,han2020suite}.

\section{Methods}

The model was deliberately kept simple so that changes in reliability could be attributed to specific binding steps. The aim was to identify patterns that may guide later system-specific studies, not to predict the behavior of one particular drug, assay, or intracellular pathway.

\subsection{Minimal ternary binding network}
\label{sec:reaction_network}

The two routes lead to the same final complex: either $A$ binds mediator first, or $B$ does. Keeping both routes allows the model to examine whether molecules accumulate as separate binary complexes rather than completing productive assembly. The model contains two binding partners, $A$ and $B$, a mediator $R$, two binary complexes, $AR$ and $BR$, and the productive ternary complex $ARB$:
\begin{equation}
S=\{A,B,R,AR,BR,ARB\}.
\end{equation}
The complex $ARB$ can form through either binary intermediate:
\begin{equation}
\begin{aligned}
A+R &\xrightleftharpoons[k_{-1}]{k_{+1}} AR, \qquad
B+R \xrightleftharpoons[k_{-2}]{k_{+2}} BR,\\
AR+B &\xrightleftharpoons[k_{-3}]{k_{+3}} ARB, \qquad
BR+A \xrightleftharpoons[k_{-4}]{k_{+4}} ARB.
\end{aligned}
\end{equation}
The model retains the two processes required to produce a hook effect: the formation of productive $ARB$ and the diversion of binding partners into separate binary complexes. The total amounts of $A$, $B$ and $R$ are conserved:
\begin{equation}
\begin{aligned}
A_{\rm tot}&=[A]+[AR]+[ARB],\\
B_{\rm tot}&=[B]+[BR]+[ARB],\\
R_{\rm tot}&=[R]+[AR]+[BR]+[ARB].
\end{aligned}
\end{equation}
These equations define the material available in each simulation. At baseline, the two routes to $ARB$ were thermodynamically balanced, with $K_{d,AR}K_{d,3}=K_{d,BR}K_{d,4}$. Because the network has a deficiency $\delta=1$, stability and convergence were checked numerically rather than assumed from a deficiency-zero theorem \cite{gunawardena2003chemical,feinberg1987chemical,feinberg2019foundations}. Further details are given in Supplementary Information (SI) Note S1. 

\subsection{Deterministic simulations}
\label{sec:deterministic_methods}

The deterministic calculation treats concentrations as continuous quantities. Its steady-state result describes the expected amount after the reaction has settled, but it does not indicate whether a particular low-copy compartment contains zero, one, or several complexes. The deterministic equations describe the concentrations expected for a well-mixed system. They are listed in SI Note S2. At each mediator dose $R_{\rm tot}$, simulations began with free $A$, $B$, and $R$ and no preformed complexes. The equations were integrated until changes in species concentrations became negligible. The final ternary-complex concentration was denoted by $[ARB]^*$ and used to draw the hook curve and compare against stochastic results.

Convergence was assessed by comparing changes in species concentrations over time. Representative trajectories reached their final values within about $120~{\rm s}$, and longer runs to $500$--$1000~{\rm s}$ were used as conservative checks. Stability was also checked numerically within the allowed concentration space set by $A_{\rm tot}$, $B_{\rm tot}$, and $R_{\rm tot}$. The limiting behavior of the hook curve was verified, where at low mediator dose, $[ARB]^*$ increases approximately with $R_{\rm tot}$, while at high mediator dose, $[ARB]^*$ decreases approximately as $R_{\rm tot}^{-1}$ because excess mediator separates $A$ and $B$ into binary complexes (SI Note S2).

\subsection{Stochastic simulations and reliability metrics}
\label{sec:stochastic_methods}

Unlike concentration equations, SSA tracks individual binding and unbinding events and therefore yields an integer number of molecules at each time point. This distinction makes it possible to ask how consistently a given dose works across repeated, otherwise identical compartments. The same network was simulated using the exact stochastic simulation algorithm developed by Gillespie (SSA) \cite{gillespie1977exact}. In femtoliter volumes, even a nanomolar concentration can correspond to only a few molecules. A nonzero average concentration therefore does not mean that every compartment contains $ARB$ \cite{bhalla2004small,hahl2016comparison}. Simulations used volumes of $1$, $10$ and $100~{\rm fL}$ (femtoliters). Deterministic solutions provided average references; individual stochastic trajectories were used to measure success and failure (SI Note S3).

To evaluate these results, endpoint extinction was defined as the probability ($P$) that no productive complex was present at time $T$:
\begin{equation}
P_{\rm ext}(T)=P\left(N_{ARB}(T)=0\right).
\end{equation}
Here, $N$ corresponds to the copy number. Endpoint success was the probability that at least one $ARB$ molecule was present:
\begin{equation}
P\left(N_{ARB}(T)\geq 1\right)=1-P_{\rm ext}(T).
\end{equation}
For each $R_{\rm tot}$, many SSA trajectories were simulated, and a trajectory was counted as successful if $N_{ARB}(T)\geq1$. The success probability was estimated as
\begin{equation}
P\left(N_{ARB}(T)\geq1\right)
=
\frac{\text{number of successful trajectories}}
{\text{total number of trajectories}}.
\end{equation}
The 90\% reliability window was defined as the tested mediator-dose range over which at least 90\% of trajectories contained one or more $ARB$ molecules at the endpoint:
\begin{equation}
P\left(N_{ARB}(T)\geq 1\right)\geq 0.90 .
\end{equation}
The main stochastic endpoint was $T=100~{\rm s}$. A dose met the reliability criterion when at least 90\% of trajectories satisfied the chosen threshold. To distinguish a single complex from a larger amount, the analysis also used $N_{ARB}(T)\geq N_{\rm crit}$ with $N_{\rm crit}=5$ or $10$. These thresholds are operational measures within the model, not experimentally established levels of biological activity.

Timing was quantified using the first-passage time to the chosen copy-number threshold:
\begin{equation}
T_{\rm crit}
=
\inf\left\{
t\geq0:
N_{ARB}(t)\geq N_{\rm crit}
\right\}.
\end{equation}
For example, if a trajectory first reached five $ARB$ molecules at $32~{\rm s}$, its first-passage time for $N_{\rm crit}=5$ was $32~{\rm s}$. If it never reached five by the endpoint, it was recorded as not having crossed the threshold. Thus, $P(T_{\rm crit}\leq t)$ measures the fraction reaching the threshold at least once by time $t$. Endpoint mediator molecules were also divided into free $R$, binary-bound $AR+BR$, and productive $ARB$ to identify where the material accumulated. The two measures are complementary: endpoint success asks what remains at the observation time, whereas first passage asks whether the required number appeared at any earlier time. A complex may reach the threshold and subsequently dissociate, so the two probabilities need not be equal.

\subsection{Design perturbations and asymmetry parameterizations}
\label{sec:design_catalytic_methods}

These changes act at different stages of assembly. The first changes how quickly the productive complex appears; the second changes how long it survives. Both are tested against the unchanged baseline to determine whether either feature limits success. Changes to binding-reaction rates were tested for their ability to improve reliable assembly. The kinetic factor $\alpha_{\rm on}$ increased the rate at which an existing binary complex completed to $ARB$. The stability factor $\alpha_{\rm off}$ reduced the rate at which formed $ARB$ broke apart. These are separate, idealized changes to the reaction rates. After either $AR$ or $BR$ had formed, the ternary-complex association rates were scaled as
\begin{equation}
k_{+3}=\alpha_{\rm on} k_{+3}^{(0)},\qquad
k_{+4}=\alpha_{\rm on} k_{+4}^{(0)} .
\end{equation}
Larger $\alpha_{\rm on}$ therefore increased ternary completion. The stability change reduced the reverse rates:
\begin{equation}
k_{-3}= \frac{k_{-3}^{(0)}}{\alpha_{\rm off}},
\qquad
k_{-4}= \frac{k_{-4}^{(0)}}{\alpha_{\rm off}} .
\end{equation}
Larger $\alpha_{\rm off}$ reduced breakup of $ARB$. The two factors were varied separately and together. The kinetic-only case had $\alpha_{\rm off}=1$; the stability-only case had $\alpha_{\rm on}=1$. The correlated combined sweep set $\alpha_{\rm on}=\alpha_{\rm off}$. One-dimensional tests used a mediator dose below the deterministic hook peak, where baseline finite-copy assembly was unreliable. Changes to the first binding step, before formation of the ternary complex, were tested separately. The asymmetry parameter $\gamma$ changed whether the $AR$ or $BR$ route was favored. This differs from $\alpha_{\rm on}$ and $\alpha_{\rm off}$, which act after binary binding. The asymmetry parameter was defined as
\begin{equation}
\gamma=\frac{K_{d,AR}}{K_{d,BR}}.
\end{equation}
At $\gamma=1$, the first-step affinities are equal. The geometric mean of $K_{d, AR}$ and $K_{d, BR}$ was held fixed, while the second-step rates were unchanged. This scan therefore tests sensitivity to route preference; away from $\gamma=1$, it generally fails to satisfy the thermodynamic route-balance condition. A balanced alternative is specified in SI Note S4 but was not used for the plotted asymmetry scan.

\subsection{Catalytic ODE model and kinetic readouts}
\label{sec:catalytic_methods}

This extension connects assembly to an observable consequence. Even if two mediator doses produce similar amounts of $ARB$ at long times, they may differ in product output if one forms the complex later during the observation period. The catalytic extension tested whether formed $ARB$ could produce an output rather than remain assembled. Product $Y$ represents the converted form of $A$ and was generated through
\begin{equation}
ARB\rightarrow B+R+Y,
\end{equation}
Product formation therefore required the presence of $ARB$. The catalytic model was simulated deterministically to examine average output; stochastic failures were evaluated in the non-catalytic model. The amount formed by $120~{\rm s}$ was $Y_{120}=Y(120~{\rm s})$. The time to convert half of the initial $A$ was
\begin{equation}
T_{Y,50}=\inf\{t:Y(t)\geq0.5A_{\rm tot}\}.
\end{equation}
The early product-formation rate was estimated from the slope of $Y(t)$ over $0$--$10~{\rm s}$.
Because product formation requires the productive ternary complex, two assembly-time readouts were used. The first readout measured how quickly each mediator dose approached its own non-catalytic steady-state level:
\begin{equation}
\tau_{\rm ARB,local}
=
\inf\{t:[ARB](t)\geq0.9[ARB]^*(R_{\rm tot})\}.
\end{equation}
Here, $[ARB]^*(R_{\rm tot})$ is the non-catalytic steady-state level at that particular dose. This first measure reports how quickly a dose reaches 90\% of its own final level, which may still be small. The second uses half of the maximum steady-state level across all doses as a common threshold:
\begin{equation}
\tau_{\rm ARB,use}
=
\inf\left\{t:[ARB](t)\geq0.5\max_{R_{\rm tot}}[ARB]^*(R_{\rm tot})\right\}.
\end{equation}
Unlike the local measure, $\tau_{\rm ARB,use}$ asks when a dose reaches the same absolute threshold. The 50\% value is an illustrative model criterion, not a validated minimum for product generation. Its role is to distinguish rapid approach to a low local level from formation of an appreciable amount of complex (SI Note S5).

\section{Results}

The results first establish the average hook effect and then examine why the same dose can give different outcomes in small compartments. The later analyses test changes in the binding steps and ask whether the assembled complex can generate product on time.

\subsection{Deterministic hook behavior becomes finite-volume stochastic failure}

The deterministic curve is independent of compartment volume when starting concentrations and rate constants are unchanged. In the stochastic simulations, however, reducing volume reduces the expected number of molecules and makes zero-copy outcomes more likely. Comparing volumes therefore tests a limitation of the average curve rather than a change in the underlying chemical reactions.

\begin{figure}[H]
\centering
\includegraphics[width=0.98\textwidth]{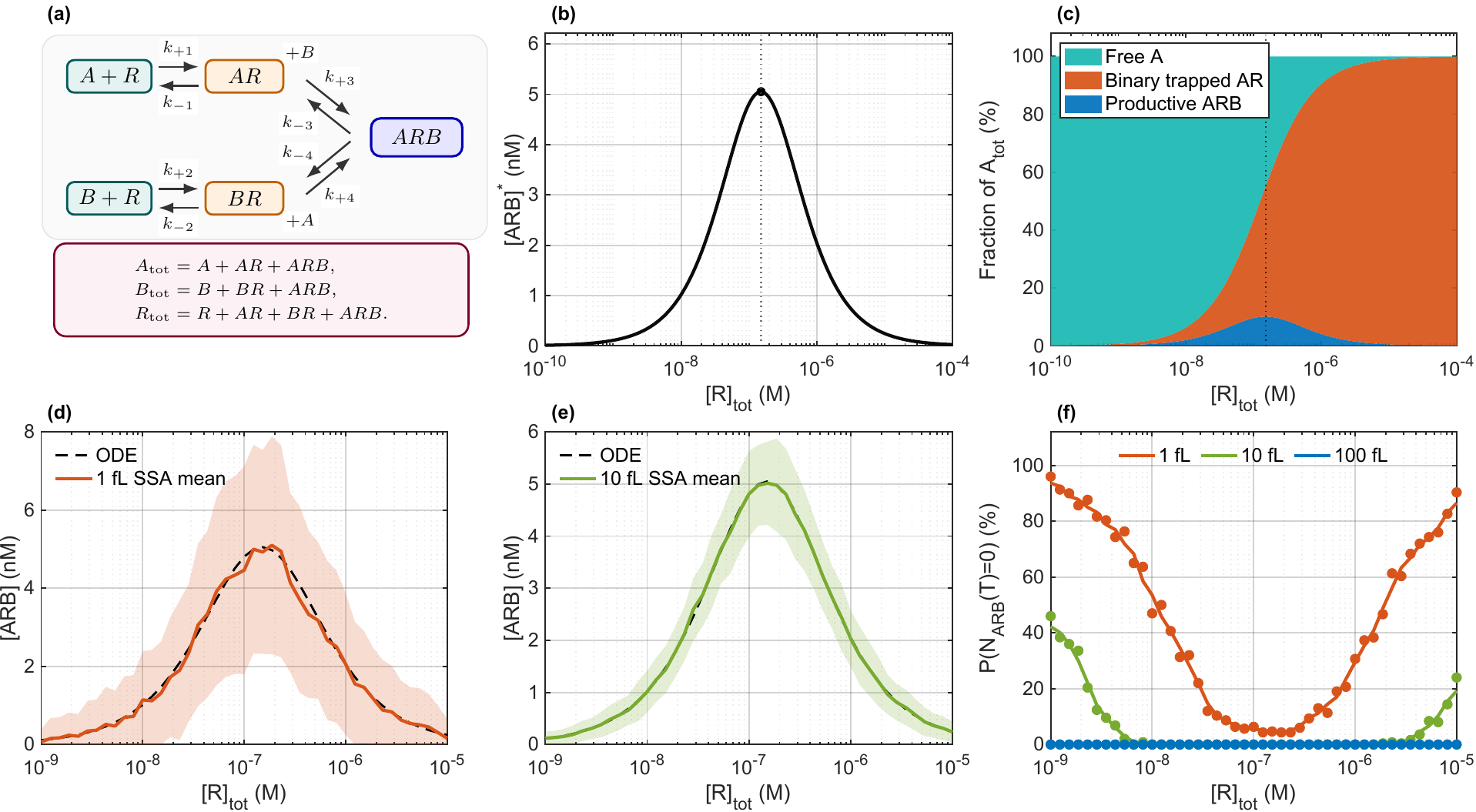}
\caption{\textbf{Deterministic hook behavior and volume-dependent stochastic reliability.}
(a) Reaction network and conservation relations. (b) Deterministic steady-state hook curve. Here, $^*$ indicates the steady-state value. (c) Partitioning of $A_{\rm tot}$. (d,e) ODE--SSA comparison at $1~{\rm fL}$ and $10~{\rm fL}$; shading shows one standard deviation. (f) Extinction probability from SSA. Simulation settings, including initial conditions and parameter values, are listed in SI Table~S1.}
\label{fig:main1}
\end{figure}

The deterministic dose response was first computed for the minimal network (Fig.~\ref{fig:main1}a). For each mediator dose, the fixed totals of $A$, $B$, and $R$ determined how much material could be free, binary-bound, or assembled as $ARB$. The steady-state response had the expected hook shape (Fig.~\ref{fig:main1}b). At low $R_{\rm tot}$, adding mediator increased $ARB$. At high $R_{\rm tot}$, excess mediator pushed material into $AR$ and $BR$, so $ARB$ decreased. SI Fig.~S1 confirms these low- and high-dose limits. Figure~\ref{fig:main1}c shows the same redistribution of $A_{\rm tot}$. As $R_{\rm tot}$ increases, free $A$ first shifts into productive $ARB$. At higher $R_{\rm tot}$, productive $ARB$ decreases and more $A$ remains in the binary-bound $AR$ state.

Stochastic simulations revealed how much individual compartments can differ from the average hook curve. Their mean response followed the deterministic trend, but the spread was large at $1~{\rm fL}$ (Fig.~\ref{fig:main1}d) and smaller at $10~{\rm fL}$ (Fig.~\ref{fig:main1}e). At $100~{\rm fL}$, the SSA mean was close to the deterministic curve (SI Fig.~S2a). At $1~{\rm fL}$, compartments often contained no $ARB$ at low or high mediator doses (Fig.~\ref{fig:main1}f). This was associated with mediator scarcity at low dose and nonproductive binding at high dose. Increasing volume reduced extinction and broadened the range with one-complex success (SI Fig.~S2b,c). The same ordering was observed at different endpoints (SI Fig.~S3).

\subsection{Mediator fate and threshold crossing distinguish stochastic failure mechanisms}

\begin{figure}[H]
\centering
\includegraphics[width=0.95\textwidth]{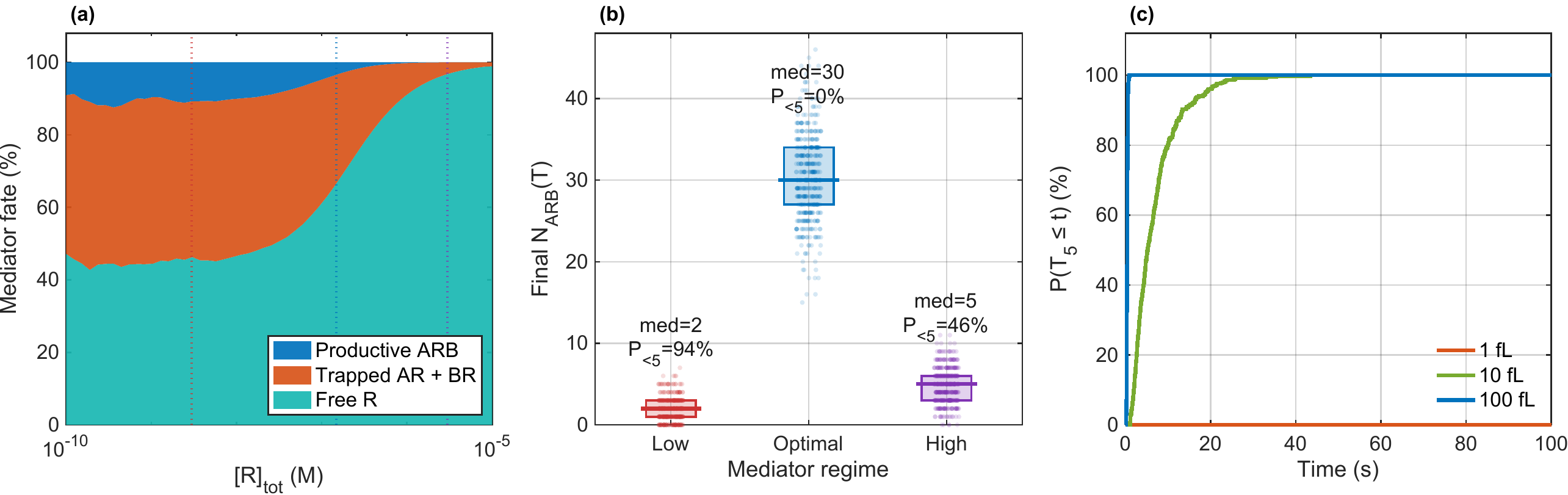}
\caption{\textbf{Stochastic failure mechanism at finite copy number.}
(a) Mediator fate map at $10~{\rm fL}$ and $T=100~{\rm s}$, showing free $R$, binary-trapped $AR+BR$, and productive $ARB$.  Vertical dotted lines mark the selected low, optimal, and high mediator regimes. 
(b) Endpoint distributions of $N_{ARB}(T)$ at $10~{\rm fL}$ with $N_{\rm crit}=5$. 
(c) First-passage probability for reaching $N_{ARB}\geq5$ at the low-mediator stress dose across $1$, $10$, and $100~{\rm fL}$. Simulation settings, including initial conditions and parameter values, are listed in SI Table~S1.}
\label{fig:main2}
\end{figure}

The same low endpoint copy number can arise for different reasons. A compartment may lack sufficient mediator to begin assembly, or it may contain abundant mediator while the partners remain separated. The species fractions and copy-number distributions help distinguish these situations. Endpoint extinction showed whether $ARB$ was absent, but not why. To distinguish the failure mechanisms, the fraction of mediator in free $R$, binary-bound $AR+BR$, or productive $ARB$ was calculated (Fig.~\ref{fig:main2}a).

Two causes of failure emerged. At low $R_{\rm tot}$, few mediator molecules were available for productive assembly. At high $R_{\rm tot}$, mediators were abundant, but most were free or in the binary complexes $AR$ and $BR$ (Fig.~\ref{fig:main2}a; SI Fig.~S4a). Endpoint copy-number distributions supported this distinction (Fig.~\ref{fig:main2}b). At the intermediate, hook-optimal dose, most trajectories contained many $ARB$ molecules. At low and high doses, many trajectories remained below five, for different reasons. Supporting distributions, threshold comparisons and example trajectories are shown in SI Fig.~S4b--d.

The next analysis examined whether trajectories ever reached a total of five $ARB$ molecules during the $100~{\rm s}$ observation period. At the low-dose stress condition, this was rare at $1~{\rm fL}$, more frequent at $10~{\rm fL}$, and nearly certain at $100~{\rm fL}$ (Fig.~\ref{fig:main2}c). Thus, failure in the smallest volume was not just a slow response: many trajectories did not reach the threshold during the observation period. The threshold and volume comparisons are shown in SI Fig.~S4e,f.

\subsection{Design perturbations reshape the finite-copy reliability window}

\begin{figure}[H]
\centering
\includegraphics[width=0.98\textwidth]{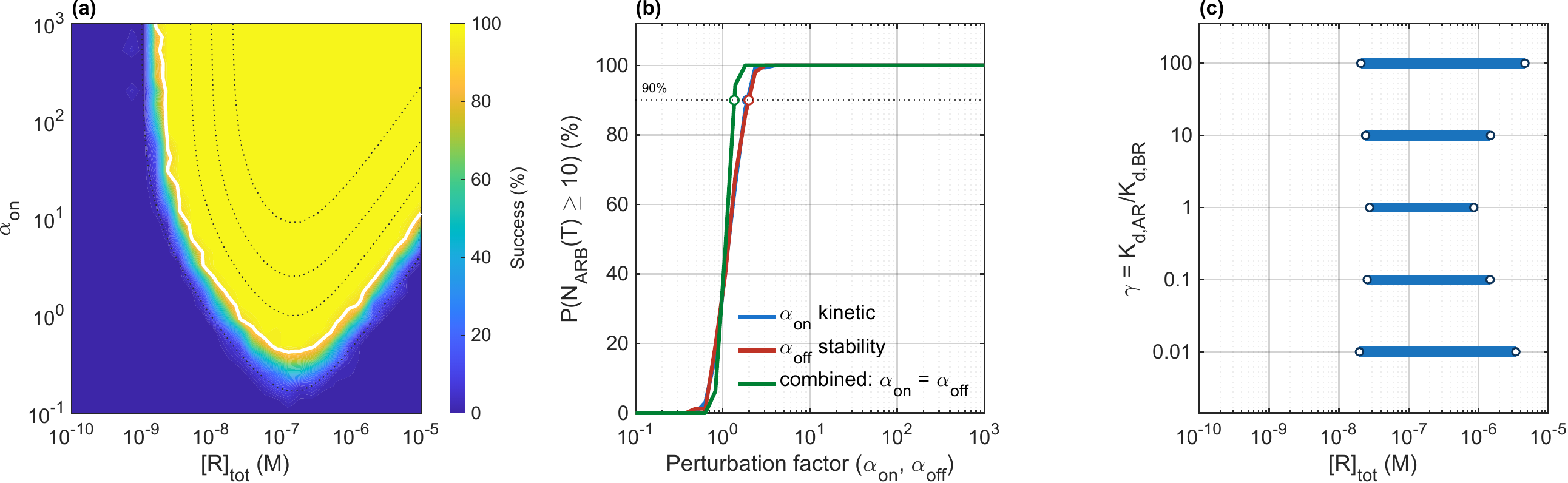}
\caption{\textbf{Design perturbations and binding asymmetry.}
(a) Endpoint success across $R_{\rm tot}$ and $\alpha_{\rm on}$ at $10~{\rm fL}$ and $T=100~{\rm s}$; the white contour marks 90\% reliability.
(b) Endpoint success for kinetic, stability, and combined perturbations at a selected low mediator dose. Markers indicate 90\% success.
(c) Reliable mediator-dose ranges for selected $\gamma$ values. Success requires $N_{ARB}(T)\geq10$. Simulation settings are listed in SI Table~S1.}
\label{fig:main3}
\end{figure}

Three perturbations were examined to improve ternary-complex assembly: faster completion into $ARB$ ($\alpha_{\rm on}$), slower dissociation of formed $ARB$ ($\alpha_{\rm off}$), and altered preference for the initial $AR$ or $BR$ binding route ($\gamma$). Increasing $\alpha_{\rm on}$ widened the mediator-dose range with at least 90\% endpoint success (Fig.~\ref{fig:main3}a). However, very low doses remained limited by mediator scarcity and high doses by nonproductive binary binding (SI Fig.~S5a--d). At the selected low-dose condition, increasing $\alpha_{\rm off}$ also improved success, as did the combined perturbation (Fig.~\ref{fig:main3}b; SI Fig.~S5e). Thus, reliability improved when $ARB$ formed faster or remained assembled longer.

First-passage analysis showed that stronger perturbations also increased the probability of reaching ten $ARB$ molecules by $100~{\rm s}$ (SI Fig.~S6). Under the selected conditions, reducing dissociation reached high success at slightly smaller perturbation factors than accelerating completion. The plotted timing summary assigns 100 s to trajectories that never reached the threshold, so its plateau does not represent an actual crossing time. Similarly, the estimated 90\% reliability boundaries should be interpreted in terms of overall trends rather than small differences between contours.

Changing $\gamma$ altered the preferred initial binding route and changed both the position and width of the reliable dose range (Fig.~\ref{fig:main3}c; SI Fig.~S7). Unlike $\alpha_{\rm on}$ and $\alpha_{\rm off}$, this perturbation did not directly modify ternary-complex completion or breakup. Because the plotted asymmetry scan did not preserve thermodynamic balance, these findings should be interpreted as sensitivity results requiring confirmation under a balanced formulation.

\subsection{The catalytic extension adds timing and product-conversion constraints}

\begin{figure}[H]
\centering
\includegraphics[width=0.98\textwidth]{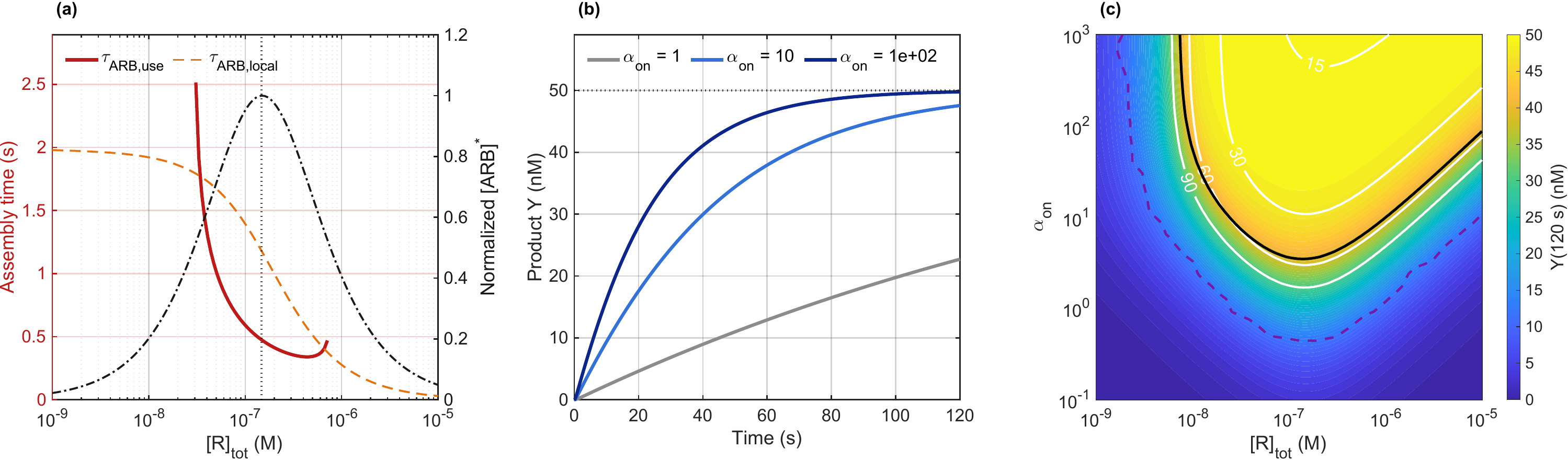}
\caption{\textbf{Kinetic and catalytic performance.}
(a) Assembly times and normalized non-catalytic steady-state $ARB$ concentration.
(b) Product accumulation at the hook-optimal mediator dose for selected $\alpha_{\rm on}$ values.
(c) Endpoint product, $Y_{120}$, across $R_{\rm tot}$ and $\alpha_{\rm on}$. White contours indicate $T_{Y,50}$, the solid black contour marks $Y_{120}=0.8A_{\rm tot}$, and the dashed magenta contour shows the non-catalytic 90\% SSA reliability boundary. Simulation settings are listed in SI Table~S1.}
\label{fig:main4}
\end{figure}

The catalytic extension examined whether formed $ARB$ could generate product $Y$ within $120~{\rm s}$. Product formation depends not only on how much $ARB$ forms, but also on when it becomes available (Fig.~\ref{fig:main4}a).

The steady-state $ARB$ concentration was highest at intermediate mediator doses (Fig.~\ref{fig:main4}a). The local assembly time, $\tau_{\rm ARB, local}$, measured how quickly each dose reached 90\% of its own steady-state level. This could be short even at high mediator doses, where little $ARB$ ultimately formed. In contrast, $\tau_{\rm ARB, use}$ used a common threshold equal to half the maximum steady-state $ARB$ concentration. This threshold was reached mainly at intermediate doses, showing that rapid assembly alone did not guarantee substantial complex formation. Increasing $\alpha_{\rm on}$ promoted faster product accumulation (Fig.~\ref{fig:main4}b) by increasing the availability of $ARB$ for conversion into $Y$ (SI Fig.~S8a). Product output remained highest at intermediate mediator doses and decreased when mediator was either scarce or excessive (Fig.~\ref{fig:main4}c; SI Fig.~S8b). Increasing $\alpha_{\rm on}$ also expanded the high-output region and shortened the time to half conversion, consistent with the early-rate and product-trajectory results (SI Fig.~S8c,d). The overlaid contours in Fig.~\ref{fig:main4}c compare two different outcomes: deterministic product conversion by $120~{\rm s}$ and stochastic $ARB$ assembly reliability at $100~{\rm s}$. They indicate how the predicted operating ranges relate to one another, rather than defining equivalent success criteria.

\section{Discussion}

The present study demonstrates that the deterministic hook response identifies mediator concentrations that favor ternary binding but does not fully capture the reliability of complex formation when molecule numbers are small. Stochastic simulations showed that, at low doses, complex formation was limited by mediator availability, whereas excess mediator favored separate binary complexes and reduced productive binding. These effects were more pronounced in smaller compartments, where the number of $ARB$ complexes varied considerably between simulations. Thus, the mediator concentration producing the highest steady-state $[ARB]^*$ does not necessarily ensure reliable formation of a specified number of complexes within the observation period.

To examine these differences, endpoint measures assessed how many complexes remained at the observation time, whereas first-passage analysis determined whether a specified number had formed by that time. The distribution of mediator among free, binary-bound, and ternary-bound states helped explain assembly failure. Together, these measures identified mediator-dose ranges meeting the chosen success criterion. These ranges generally widened with increasing compartment volume but narrowed when more complexes were required. Because the same copy-number thresholds were used across volumes, they represented different concentration requirements.

The analysis then examined whether changes in binding rates could improve reliability. Increasing $\alpha_{\rm on}$ accelerated $ARB$ formation and widened the reliable dose range, whereas increasing $\alpha_{\rm off}$ slowed complex breakup and improved success at the selected low-dose condition. Both changes reduced the effective second-step dissociation constants but affected complex formation and persistence in different ways. Changing the relative first-step binding affinities also altered the position and width of the reliable dose range. However, this analysis did not preserve thermodynamic balance between the two routes to $ARB$ when $\gamma\neq1$. These findings therefore require confirmation using a balanced model before being applied to physical binding systems. The catalytic extension showed that product formation depends on both the amount and timing of $ARB$ formation. Faster ternary completion increased product accumulation over $120~{\rm s}$ and shortened the time to half-conversion under the tested conditions. The non-catalytic assembly-time measures further showed that rapidly approaching a low steady-state concentration does not necessarily indicate substantial complex formation. Because the catalytic calculations were deterministic, they do not establish the reliability of product formation in individual low-copy compartments.

The model assumes a well-mixed system with reactions following mass-action kinetics. It does not explicitly include spatial transport, molecular crowding, competing binding partners, or differences between cells, all of which may influence complex formation \cite{puskar2007spatial,miklos2013competitive,stoeger2016passive}. The changes in binding rates and the selected success thresholds were also used to compare model conditions rather than represent experimentally established requirements. In addition, the estimated reliability boundaries are uncertain because they were calculated from a finite number of stochastic simulations. Applying these findings to specific biochemical systems will require measured reaction rates, thermodynamically consistent binding parameters, and success criteria based on the intended biological function. Extending the catalytic model to stochastic simulations would also allow the reliability of product formation to be examined directly.

\section{Conclusions}

The present study shows that a favorable average ternary-binding response does not necessarily guarantee reliable complex formation when the number of molecules is small. Stochastic simulations identified how mediator concentration, compartment volume, the required number of complexes, and binding rates affect assembly reliability. Faster ternary formation and slower complex breakup improved reliability under the tested conditions, while the catalytic analysis showed that the timing of complex formation also influences product output. These findings provide a basis for evaluating both the extent and reliability of ternary complex formation, although their application to specific biological systems requires experimentally measured parameters and appropriate success criteria.

\section*{Acknowledgements}

The author thanks Eduardo Sontag for helpful discussions.

\section*{Funding}

This work was supported by the Mathematical Research Impact Centric Support (MATRICS) grant from ANRF (formerly SERB), India (MTR/2022/000369), and the IIT Bombay Seed Grant (RD/0523-IRCCSH0-024).

\section*{Ethics}

This study used only mathematical modeling and simulation. It did not involve human participants, human tissue, animals, or fieldwork.

\section*{Data accessibility}

This study generated no experimental data. 

\section*{Competing interests}

The author declares no competing interests.

\section*{Author Contributions}

Deepak K. Agrawal: Conceptualization, methodology, software, formal analysis, validation, visualization, writing--original draft, writing--review and editing, funding acquisition, and project administration.

\section*{Supplementary material}

Supplementary Information includes derivations, implementation details, parameter tables, and additional figures.

\section*{Declaration of generative AI use}
ChatGPT (OpenAI) was used to assist with language editing, formatting, reference checking, and code review. The author reviewed the suggestions and takes full responsibility for the scientific content and conclusions.

\clearpage
\begingroup
\setstretch{1.02}
\setcounter{equation}{0}
\renewcommand{\theequation}{S\arabic{equation}}
\setcounter{figure}{0}
\renewcommand{\thefigure}{S\arabic{figure}}
\setcounter{table}{0}
\renewcommand{\thetable}{S\arabic{table}}
\pdfbookmark[1]{Supplementary Information}{supplementary-info}
\begin{center}
\vspace*{1em}
{\LARGE\bfseries Supplementary Information\par}
\vspace{1.25em}
{\large From the Deterministic Hook Effect to Stochastic Reliability in Ternary Complex Formation\par}
\vspace{1.25em}
{\normalsize Deepak K. Agrawal\par}
\vspace{0.5em}
{\normalsize Department of Biosciences and Bioengineering,\par}
{\normalsize Indian Institute of Technology Bombay, Mumbai, India\par}
\end{center}
\clearpage

\section*{Supplementary Note S1. Network structure and thermodynamic consistency}

The reversible six-species network conserves the total amounts of $A$, $B$ and $R$. In copy-number form, these balances are
\begin{equation}
\begin{aligned}
N_{A,{\rm tot}}&=N_A+N_{AR}+N_{ARB},\\
N_{B,{\rm tot}}&=N_B+N_{BR}+N_{ARB},\\
N_{R,{\rm tot}}&=N_R+N_{AR}+N_{BR}+N_{ARB}.
\end{aligned}
\end{equation}
Conservation does not by itself guarantee thermodynamic consistency between the two binding routes. Defining
\begin{equation}
K_{d,AR}=\frac{k_{-1}}{k_{+1}},\quad
K_{d,BR}=\frac{k_{-2}}{k_{+2}},\quad
K_{d,3}=\frac{k_{-3}}{k_{+3}},\quad
K_{d,4}=\frac{k_{-4}}{k_{+4}},
\end{equation}
The route-balance condition is
\begin{equation}
K_{d,AR}K_{d,3}=K_{d,BR}K_{d,4}.
\end{equation}
The symmetric baseline and the second-step design perturbations satisfy this condition. The first-step asymmetry scan, which keeps second-step rates fixed, generally does not (SI Note S4). The seven distinct reaction complexes are
\begin{equation}
\mathcal C=\{A+R,AR,B+R,BR,AR+B,BR+A,ARB\}.
\end{equation}
They form three linkage classes:
\begin{equation}
\begin{aligned}
\mathcal L_1 &: A+R\leftrightarrow AR,\\
\mathcal L_2 &: B+R\leftrightarrow BR,\\
\mathcal L_3 &: AR+B\leftrightarrow ARB\leftrightarrow BR+A.
\end{aligned}
\end{equation}
The three independent conservation laws give a stoichiometric rank of $r=6-3=3$. Hence, with $n_c=7$ and $\ell=3$, the deficiency is
\begin{equation}
\delta=n_c-\ell-r=7-3-3=1.
\end{equation}
The deficiency-zero theorem is therefore not applicable; convergence was assessed numerically.

\clearpage
\section*{Supplementary Note S2. Deterministic model and hook-effect limits}

The mass-action equations are
\begin{equation}
\begin{aligned}
\frac{d[A]}{dt}&=-k_{+1}[A][R]+k_{-1}[AR]-k_{+4}[BR][A]+k_{-4}[ARB],\\
\frac{d[B]}{dt}&=-k_{+2}[B][R]+k_{-2}[BR]-k_{+3}[AR][B]+k_{-3}[ARB],\\
\frac{d[R]}{dt}&=-k_{+1}[A][R]+k_{-1}[AR]-k_{+2}[B][R]+k_{-2}[BR],\\
\frac{d[AR]}{dt}&=k_{+1}[A][R]-k_{-1}[AR]-k_{+3}[AR][B]+k_{-3}[ARB],\\
\frac{d[BR]}{dt}&=k_{+2}[B][R]-k_{-2}[BR]-k_{+4}[BR][A]+k_{-4}[ARB],\\
\frac{d[ARB]}{dt}&=k_{+3}[AR][B]+k_{+4}[BR][A]-(k_{-3}+k_{-4})[ARB].
\end{aligned}
\end{equation}
Simulations began with free $A$, $B$, and $R$, with no bound complexes. To illustrate the hook-effect limits, take $A_{\rm tot}=B_{\rm tot}=C$, $K_{d,AR}=K_{d,BR}=K_d$ and $K_{d,3}=K_{d,4}=K_{d,{\rm tern}}$. At low mediator dose, $[A]\approx[B]\approx C$, giving
\begin{equation}
[ARB]^*\approx \frac{C^2[R]}{K_dK_{d,{\rm tern}}}\propto R_{\rm tot}.
\end{equation}
At high mediator dose, $[R]\approx R_{\rm tot}$ and $[AR]\approx[BR]\approx C$, so
\begin{equation}
[ARB]^*\approx\frac{C^2K_d}{K_{d,{\rm tern}}R_{\rm tot}}\propto R_{\rm tot}^{-1}.
\end{equation}
These limiting expressions explain the slopes in SI Fig.~S1. All plotted values were obtained using the full equations.

\clearpage
\section*{Supplementary Note S3. Stochastic simulation and reliability measures}

Molecule counts were related to concentrations by $N_i=[i]VN_{\rm Av}$. Initial counts were rounded, except at low mediator inputs where initial mediator counts were Poisson sampled (Table~\ref{tab:SI_parameters_and_settings}). Exact stochastic simulation algorithm (SSA) trajectories used propensities such as
\begin{equation}
a_{A+R\rightarrow AR}=\frac{k_{+1}N_AN_R}{VN_{\rm Av}},\qquad
a_{AR\rightarrow A+R}=k_{-1}N_{AR}.
\end{equation}
Endpoint success and the logarithmic width of the tested reliable dose interval were defined as
\begin{equation}
P_{\rm th}(T)=P\!\left(N_{ARB}(T)\geq N_{\rm crit}\right),\qquad
W=\log_{10}\!\left(\frac{R_{\max}}{R_{\min}}\right).
\end{equation}
The reliable interval was identified on the tested dose grid. Its estimated boundaries are not confidence limits. The percentage difference in SI Fig.~S2 was calculated as
\begin{equation}
{\rm Bias}=100\,\frac{\langle[ARB](T)\rangle_{\rm SSA}-[ARB]^*_{\rm ODE}}{[ARB]^*_{\rm ODE}}.
\end{equation}
This compares a stochastic endpoint with a deterministic steady state so that incomplete relaxation may contribute to the difference.
\clearpage
\section*{Supplementary Note S4. Binding perturbations and asymmetry}

Mediator fractions at the endpoint were defined as
\begin{equation}
(f_{\rm free},f_{\rm binary},f_{\rm productive})=
\frac{(N_R,\,N_{AR}+N_{BR},\,N_{ARB})}{N_{R,{\rm tot}}},
\end{equation}
for trajectories containing at least one mediator molecule. The kinetic factor $\alpha_{\rm on}$ increased ternary completion, whereas $\alpha_{\rm off}$ reduced complex breakup:
\begin{equation}
\begin{aligned}
k_{+j}&=\alpha_{\rm on}k_{+j}^{(0)},&
k_{-j}&=\frac{k_{-j}^{(0)}}{\alpha_{\rm off}},\qquad j\in\{3,4\}.
\end{aligned}
\end{equation}
Single-factor tests set the other factor to one; combined tests give
\begin{equation}
K_{d,j}^{\rm eff}=\frac{K_{d,j}^{(0)}}{\alpha_{\rm on}\alpha_{\rm off}},
\qquad j\in\{3,4\}.
\end{equation}
The combined one-dimensional sweep used $\alpha_{\rm on}=\alpha_{\rm off}$; the two-dimensional sweep varied them independently. First-passage time was
\begin{equation}
T_{\rm crit}=\inf\{t\geq0:N_{ARB}(t)\geq N_{\rm crit}\}.
\end{equation}
For the plotted median in SI Fig.~S6a, non-crossing trajectories were assigned the observation endpoint. First-step asymmetry was varied using
\begin{equation}
\gamma=\frac{K_{d,AR}}{K_{d,BR}},\qquad
K_{d,AR}=K_d\sqrt{\gamma},\qquad
K_{d,BR}=\frac{K_d}{\sqrt{\gamma}}.
\end{equation}
The plotted scan kept the second-step rates fixed. A thermodynamically balanced alternative would use
\begin{equation}
K_{d,3}=\frac{K_{d,{\rm tern}}}{\sqrt{\gamma}},\qquad
K_{d,4}=K_{d,{\rm tern}}\sqrt{\gamma},
\end{equation}
but this alternative was not used to generate the asymmetry figures.

\clearpage
\section*{Supplementary Note S5. Catalytic model and assembly times}

The deterministic catalytic extension used $ARB\rightarrow B+R+Y$ with flux $v_{\rm cat}=k_{\rm cat}[ARB]$. Catalysis removes $ARB$ and produces equal amounts of $B$, $R$ and $Y$:
\begin{equation}
\left.\frac{d[ARB]}{dt}\right|_{\rm cat}=-v_{\rm cat},\qquad
\left.\frac{d[B]}{dt}\right|_{\rm cat}=
\left.\frac{d[R]}{dt}\right|_{\rm cat}=\frac{d[Y]}{dt}=v_{\rm cat}.
\end{equation}
The target balance is $A_{\rm tot}=[A]+[AR]+[ARB]+[Y]$. Product readouts were
\begin{equation}
Y_{120}=Y(120~{\rm s}),\qquad
T_{Y,50}=\inf\{t\geq0:Y(t)\geq0.5A_{\rm tot}\}.
\end{equation}
Assembly times were calculated using the non-catalytic model:
\begin{equation}
\begin{aligned}
\tau_{\rm ARB,local}&=\inf\{t\geq0:[ARB](t)\geq0.9[ARB]^*(R_{\rm tot})\},\\
[ARB]_{\rm use}&=0.5\max_{R_{\rm tot}}[ARB]^*(R_{\rm tot}),\\
\tau_{\rm ARB,use}&=\inf\{t\geq0:[ARB](t)\geq[ARB]_{\rm use}\}.
\end{aligned}
\end{equation}
These are comparison thresholds, not measured requirements for catalysis. The stochastic analyses did not include product formation.

\clearpage
\section*{Supplementary parameter table}

\begin{table}[H]
\centering
\caption{Model parameters and simulation settings.}
\label{tab:SI_parameters_and_settings}
\footnotesize
\setlength{\tabcolsep}{3.5pt}
\renewcommand{\arraystretch}{1.05}
\begin{tabular}{p{0.27\textwidth}p{0.37\textwidth}p{0.28\textwidth}}
\toprule
Item & Value or setting & Use \\
\midrule
\multicolumn{3}{l}{\textbf{Shared parameters}} \\
\midrule
Species & $A$, $B$, $R$, $AR$, $BR$, $ARB$ & Minimal reversible ternary network. \\
Conserved totals & $A_{\rm tot}=B_{\rm tot}=50~{\rm nM}$; $R_{\rm tot}$ varied logarithmically & Fixed partners; swept mediator dose. \\
Baseline association rates & $k_{+1}=k_{+2}=k_{+3}=k_{+4}=10^7~{\rm M}^{-1}{\rm s}^{-1}$ & Binding rates. \\
Baseline dissociation rates & $k_{-1}=k_{-2}=k_{-3}=k_{-4}=1~{\rm s}^{-1}$ & Baseline $K_d=100~{\rm nM}$. \\
Deterministic output & $[ARB]^*$ & Steady-state ternary concentration. \\
SSA propensities & $a_{ij}=k_+N_iN_j/(VN_{\rm Av})$; $a_i=k_-N_i$ & Discrete reaction simulation. \\
Default SSA endpoint & $T=100~{\rm s}$ & Non-catalytic reliability. \\
Copy thresholds & $N_{\rm crit}=1$, $5$, $10$ & One-complex, failure and design analyses. \\
Second-step perturbations & $k_{+3,+4}\propto\alpha_{\rm on}$; $k_{-3,-4}\propto1/\alpha_{\rm off}$ & Completion and stability (SI Note S4). \\
Catalytic reaction & $ARB\rightarrow B+R+Y$; $k_{\rm cat}=0.05~{\rm s}^{-1}$ & Deterministic product formation. \\
Catalytic readouts & $Y_{120}$, $T_{Y,50}$, $\tau_{\rm ARB,local}$, $\tau_{\rm ARB,use}$ & Product and assembly timing. \\
\midrule
\multicolumn{3}{l}{\textbf{Figure-specific settings}} \\
\midrule
Fig.~1; SI Figs.~S1--S3 & $R_{\rm tot}$ sweep; $1$, $10$, $100~{\rm fL}$ & Hook curve and one-complex reliability. \\
Fig.~2; SI Fig.~S4 & Selected low, optimal and high mediator regimes; $N_{\rm crit}=1$, $5$, $10$ & Failure mechanisms and first passage. \\
Fig.~3a,b; SI Fig.~S5 & $R_{\rm tot}$, $\alpha_{\rm on}$, $\alpha_{\rm off}$; $10~{\rm fL}$; $100~{\rm s}$; $N_{\rm crit}=10$; 80 SSA runs per map point & Design perturbations and dose windows. \\
Fig.~3b; SI Fig.~S6 & Selected low-dose stress condition; $\alpha_{\rm on}$ and $\alpha_{\rm off}$ & Threshold-crossing timing. \\
Fig.~3c; SI Fig.~S7 & $R_{\rm tot}$, $\gamma$; $10~{\rm fL}$; $100~{\rm s}$ & Unbalanced asymmetry sensitivity scan. \\
Fig.~4; SI Fig.~S8 & $R_{\rm tot}$, $\alpha_{\rm on}$; $120~{\rm s}$ & Catalytic product and assembly timing. \\
\bottomrule
\end{tabular}
\end{table}
\clearpage

\section*{Supplementary figures}

\begin{figure}[H]
\centering
\includegraphics[width=0.5\textwidth]{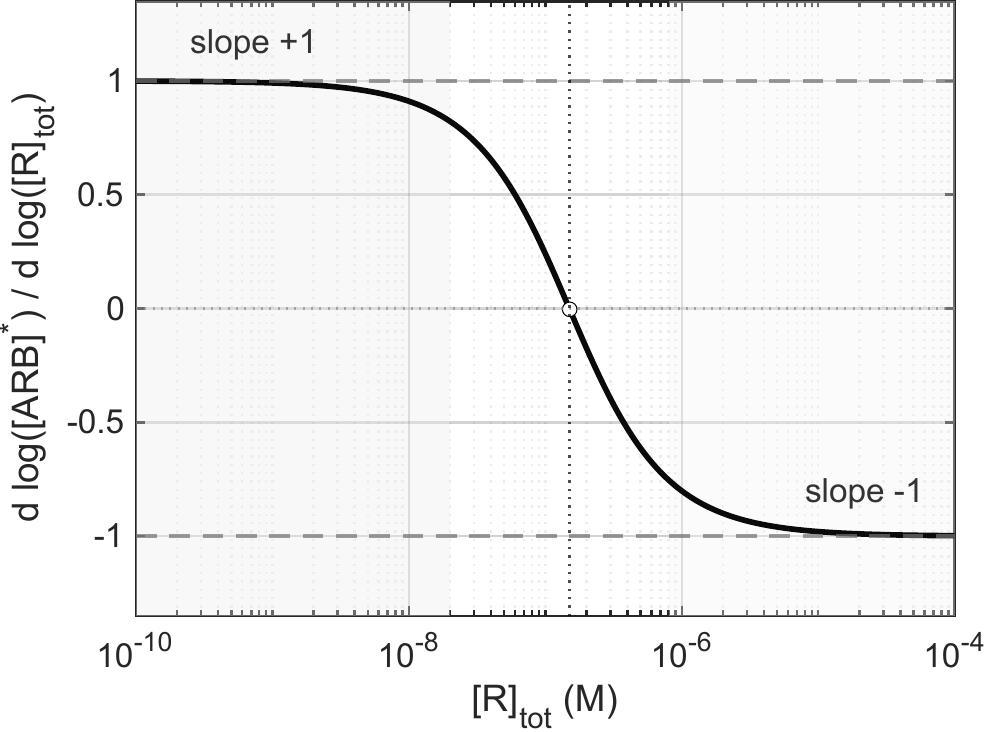}
\caption{Local slope of the deterministic hook curve. Reference slopes $+1$ and $-1$ indicate the expected low- and high-mediator limits. Baseline settings are in Table~\ref{tab:SI_parameters_and_settings}.}
\label{fig:S1}
\end{figure}

\vspace{100pt}

\begin{figure}[H]
\centering
\includegraphics[width=\textwidth]{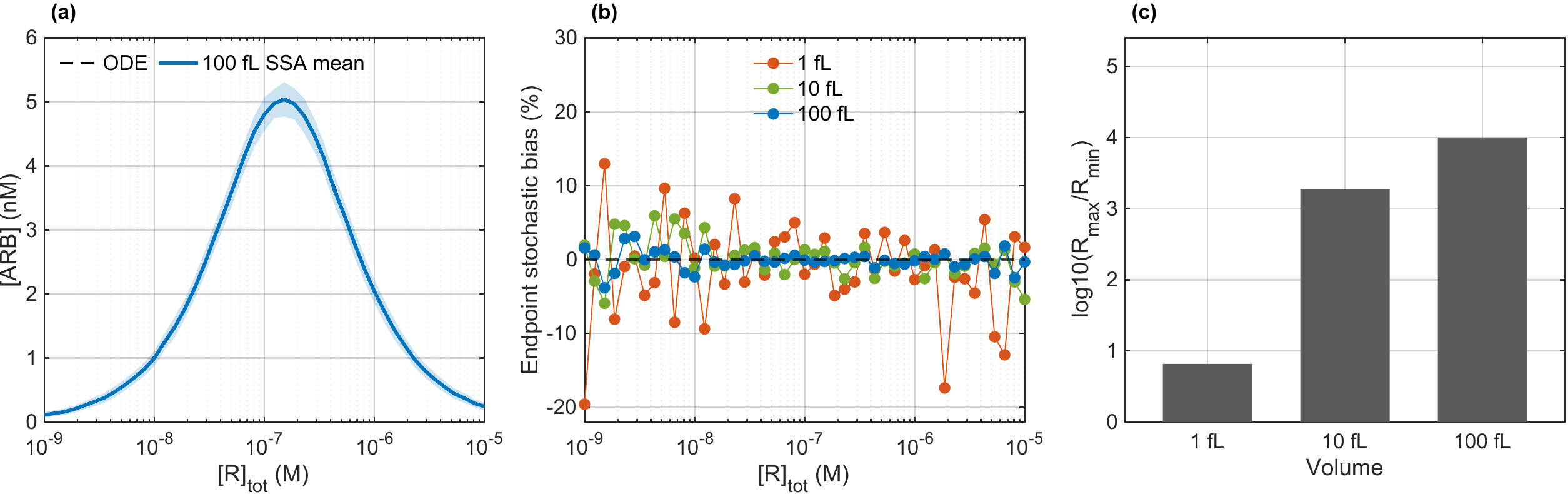}
\caption{SSA volume comparison. (a) ODE--SSA comparison at $100~{\rm fL}$, using the SSA ensemble mean at $T=100~{\rm s}$. (b) Percentage difference between the SSA endpoint mean and the ODE steady state at $1$, $10$ and $100~{\rm fL}$. (c) One-complex 90\% reliability-window width, expressed as $\log_{10}(R_{\max}/R_{\min})$.}
\label{fig:S2}
\end{figure}

\begin{figure}[H]
\centering
\includegraphics[width=\textwidth]{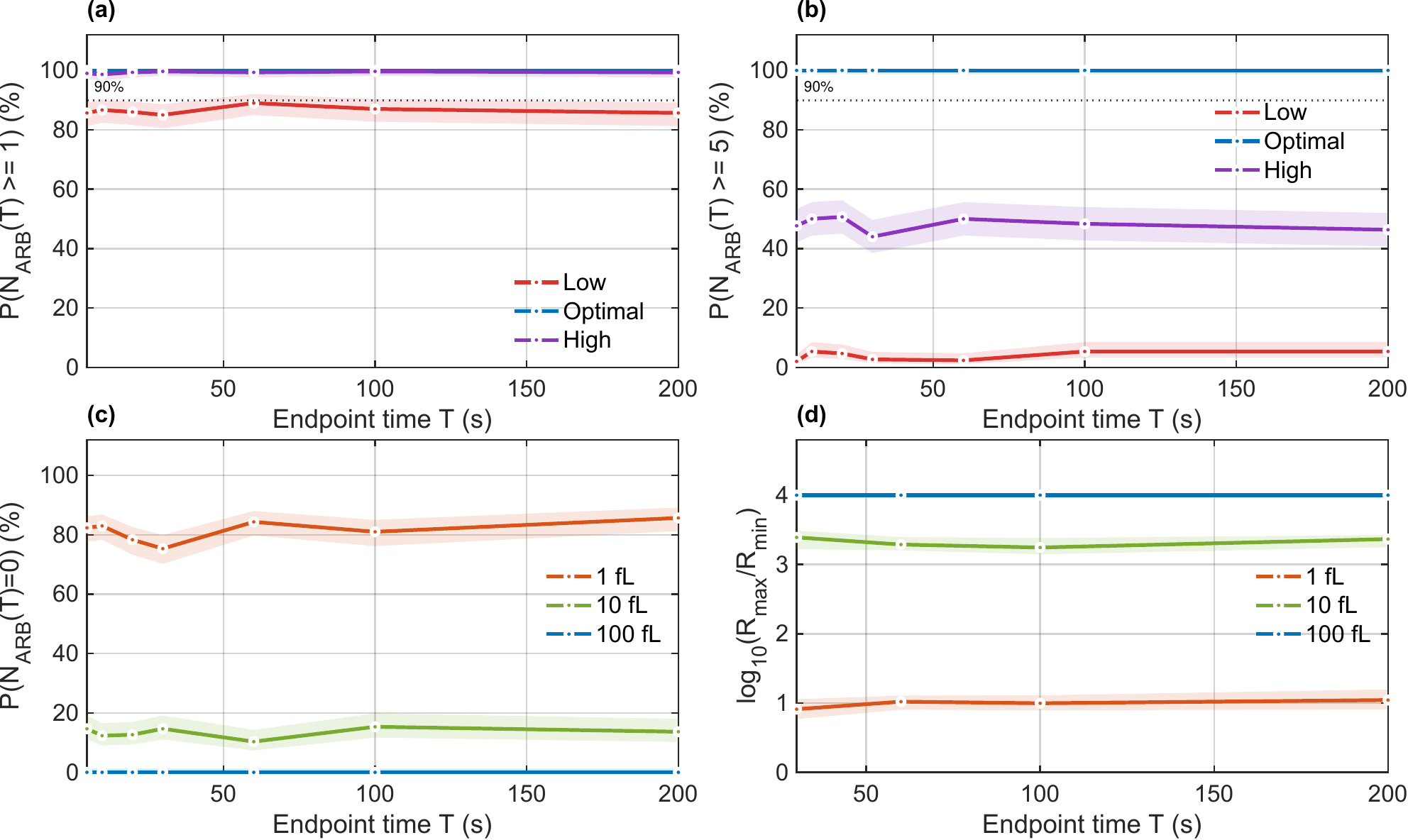}
\caption{Endpoint-time sensitivity. (a,b) Endpoint probabilities at selected low, optimal, and high mediator doses for $N_{ARB}(T)\geq1$ and $N_{ARB}(T)\geq5$. (c) Extinction at low mediator dose across $1$, $10$ and $100~{\rm fL}$. (d) One-complex 90\% reliability-window width versus endpoint time. Endpoint comparisons reuse the same SSA trajectories.}
\label{fig:S3}
\end{figure}

\vspace{40pt}
\begin{figure}[H]
\centering
\includegraphics[width=\textwidth]{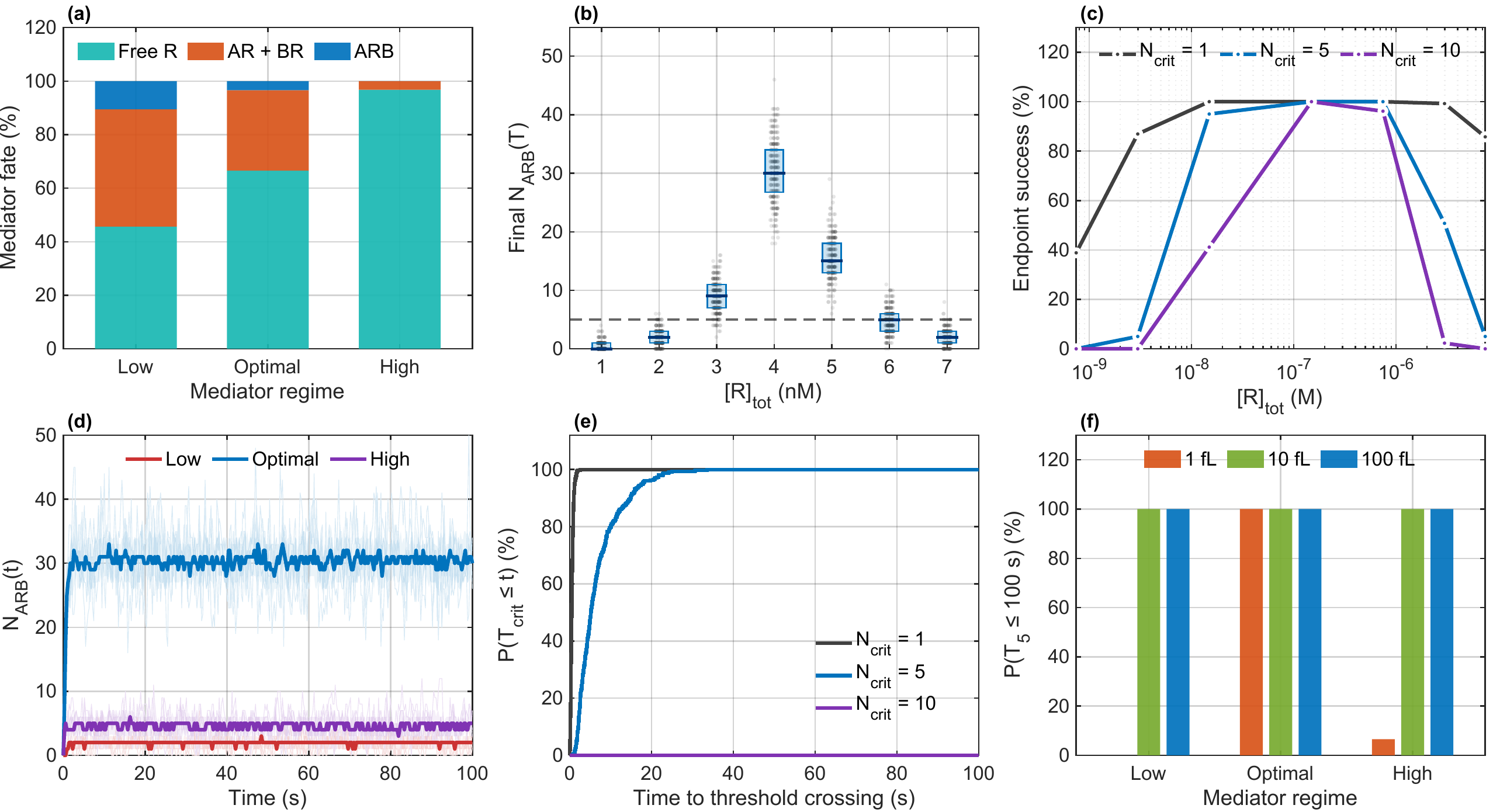}
\caption{Additional stochastic failure analyses for main Fig.~2. (a) Endpoint mediator fate at $T=100~{\rm s}$. (b) $N_{ARB}(T)$ distributions for $N_{\rm crit}=5$. (c) Endpoint threshold sensitivity at $N_{\rm crit}=1$, $5$, $10$. (d) Representative SSA trajectories. (e) First-passage curves for the same thresholds. (f) Probability of reaching $N_{ARB}\geq5$ by $100~{\rm s}$ at $1$, $10$ and $100~{\rm fL}$.}
\label{fig:S4}
\end{figure}

\begin{figure}[H]
\centering
\includegraphics[width=\textwidth]{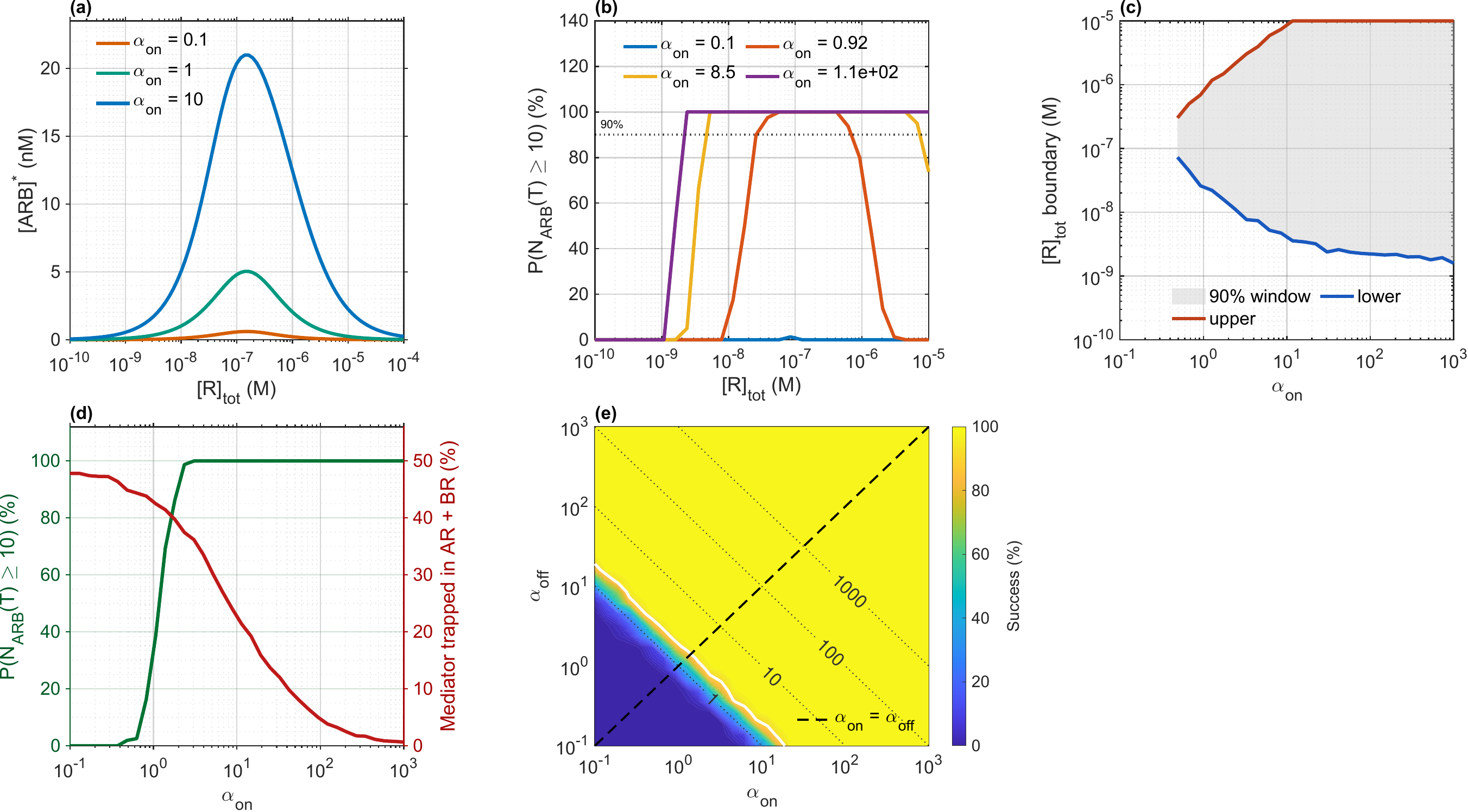}
\caption{Kinetic and combined perturbations for main Fig.~3a,b. (a) Deterministic $[ARB]^*$ for selected $\alpha_{\rm on}$. (b) Endpoint success at $10~{\rm fL}$ and $100~{\rm s}$ for $N_{ARB}(T)\geq10$. (c) 90\% reliability-window boundaries versus $\alpha_{\rm on}$. (d) Endpoint success and binary trapping at low dose. (e) Joint $\alpha_{\rm on}$--$\alpha_{\rm off}$ map with a diagonal $\alpha_{\rm on}=\alpha_{\rm off}$ and guides for constant $\alpha_{\rm on}\alpha_{\rm off}$.}
\label{fig:S5}
\end{figure}

\vspace{50pt}
\begin{figure}[H]
\centering
\includegraphics[width=\textwidth]{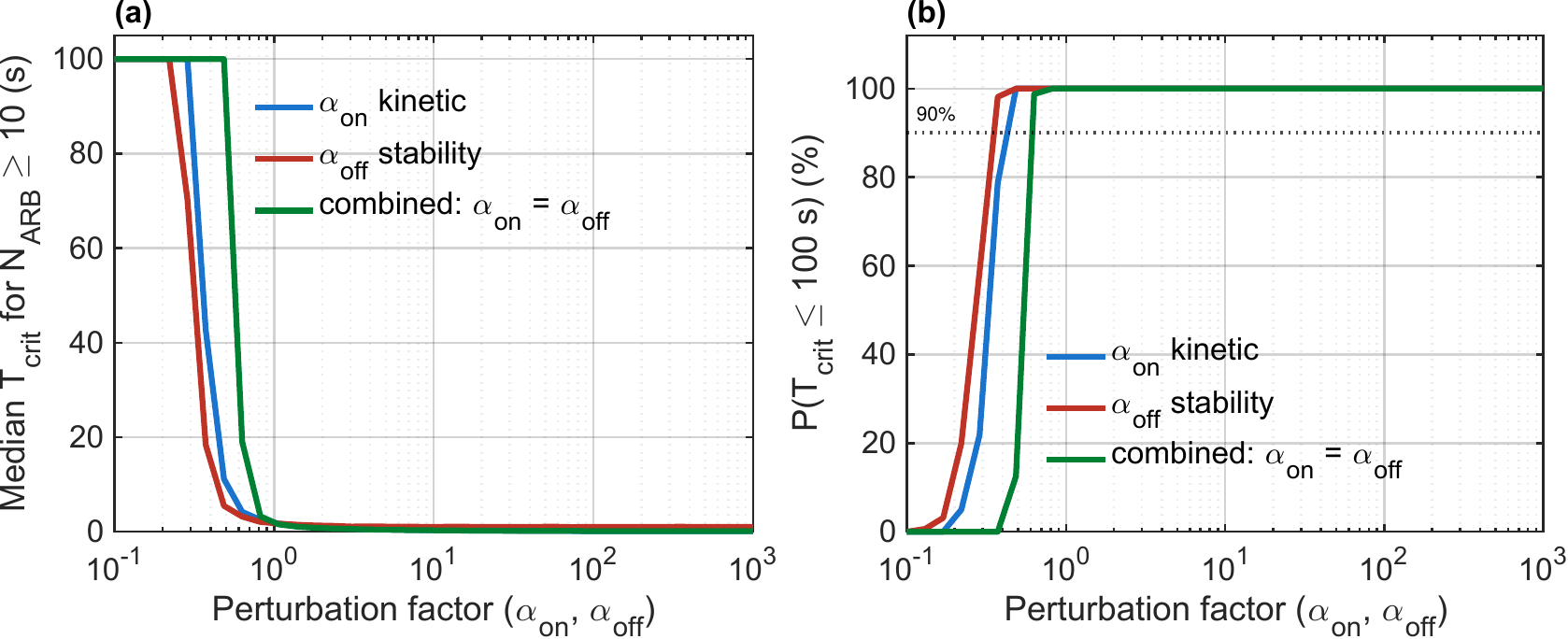}
\caption{First-passage analysis for design perturbations. (a) Median of threshold times truncated at $100~{\rm s}$; trajectories without crossing were assigned $100~{\rm s}$. (b) Probability of reaching $N_{ARB}\geq10$ by $100~{\rm s}$ at the low-dose stress condition. Kinetic, stability and combined cases vary $\alpha_{\rm on}$, $\alpha_{\rm off}$ and $\alpha_{\rm on}=\alpha_{\rm off}$, respectively.}
\label{fig:S6}
\end{figure}

\begin{figure}[H]
\centering
\includegraphics[width=\textwidth]{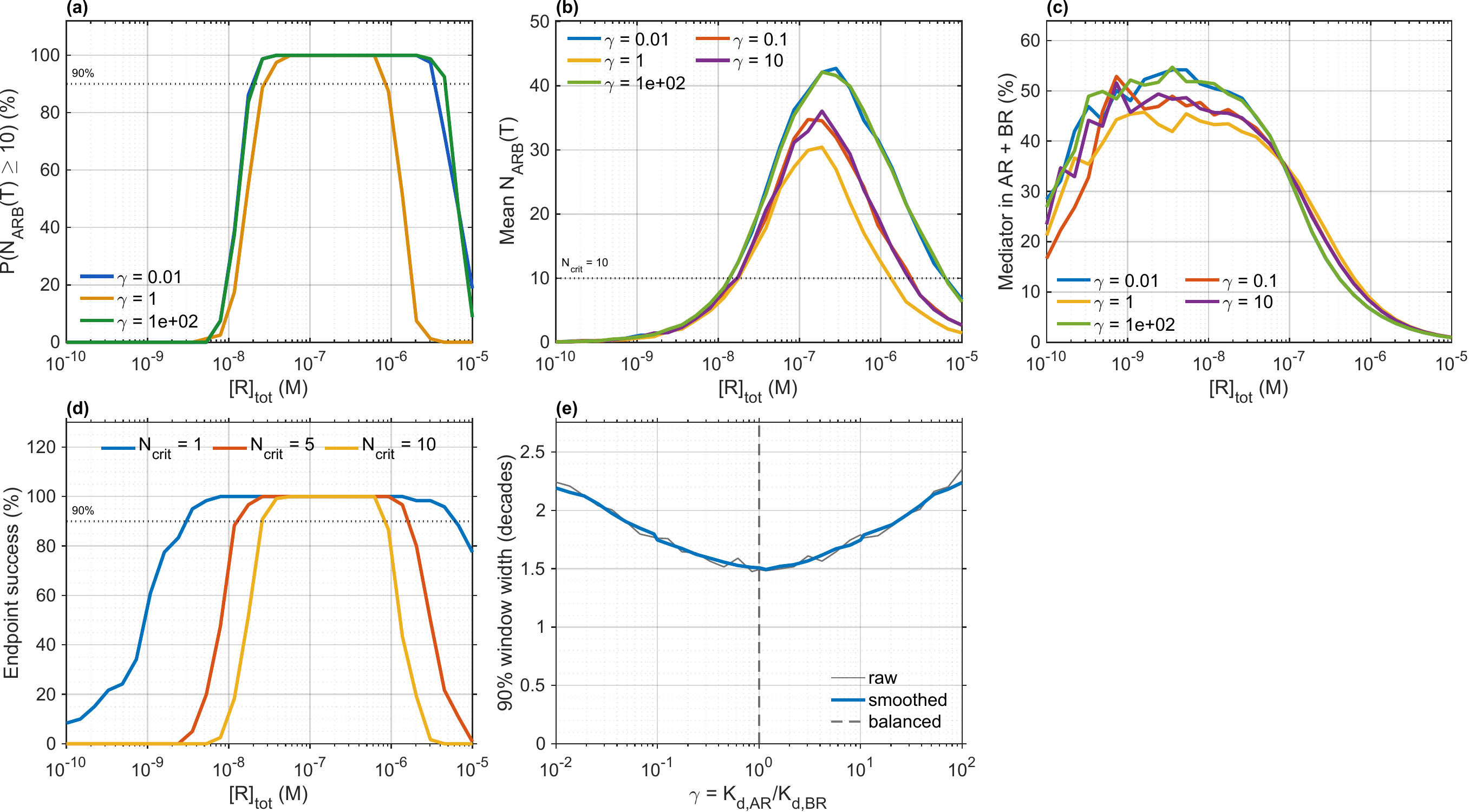}
\caption{First-step asymmetry sensitivity for main Fig.~3c. (a) Endpoint success for selected $\gamma=K_{d,AR}/K_{d,BR}$. (b) Mean $N_{ARB}(T)$, with $N_{\rm crit}=10$ as a reference. (c) Endpoint binary-trapped mediator fraction. (d) Threshold sensitivity at $\gamma=1$ for $N_{\rm crit}=1$, $5$, $10$. (e) Estimated 90\% reliability-window width versus $\gamma$ for $N_{\rm crit}=10$. The plotted parameterization does not enforce thermodynamic route balance for $\gamma\neq1$.}
\label{fig:S7}
\end{figure}

\vspace{40pt}
\begin{figure}[H]
\centering
\includegraphics[width=\textwidth]{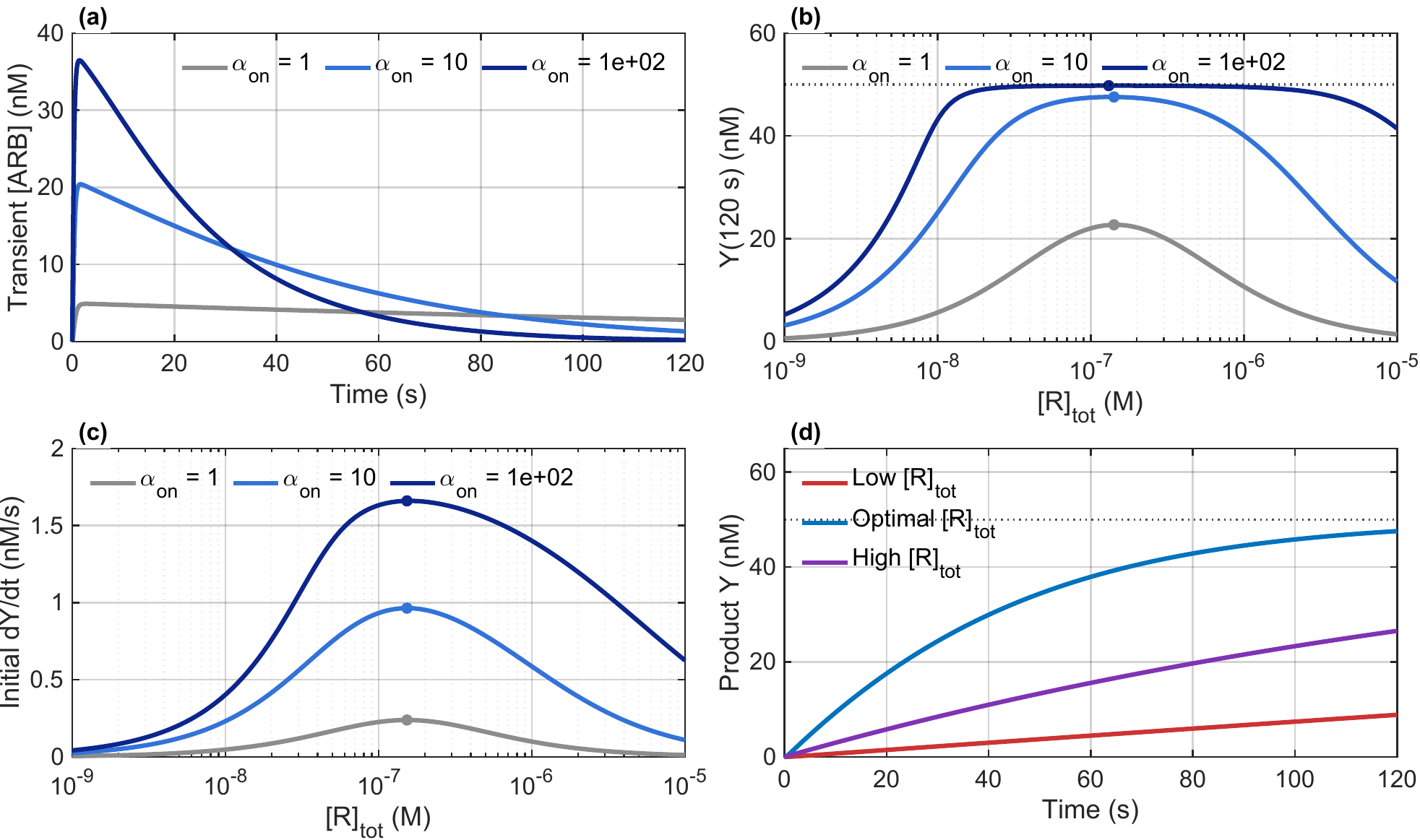}
\caption{Additional catalytic analyses for main Fig.~4. (a) Transient $[ARB](t)$ at the hook-optimal dose for selected $\alpha_{\rm on}$. (b) Endpoint product $Y_{120}=Y(120~{\rm s})$ versus $R_{\rm tot}$. (c) Initial product-formation rate estimated over $0$--$10~{\rm s}$. (d) Product trajectories at low, hook-optimal, and high mediator doses.}
\label{fig:S8}
\end{figure}
\endgroup
\end{document}